\documentclass{pspum-l}

\usepackage{amssymb}
\usepackage{cite}

\usepackage{hyperref}
\hypersetup{pdftitle={Computing Cohomology on Toric Varieties}, 
      pdfauthor={Benjamin Jurke}, 
      pdfsubject={Cohomology of line bundle; Toric variety; Alexander duality},
      pdfstartview={FitH}, pdfpagelayout={TwoColumnRight},
      bookmarksopen, bookmarksnumbered, bookmarksopenlevel=2}

\newcommand{\cohomCalg}{\href{http://wwwth.mppmu.mpg.de/members/blumenha/cohomcalg/}{\text{\fontfamily{put}\bfseries\footnotesize\selectfont cohomCalg}}}

\newcommand{\beq}{\begin{equation}}  \newcommand{\eeq}{\end{equation}}
\newcommand{\bal}{\begin{aligned}}   \newcommand{\eal}{\end{aligned}}

\def\IC{\mathbb{C}}
\def\IZ{\mathbb{Z}}
\def\IN{\mathbb{N}}

\def\cN{\mathcal{N}}
\def\cO{\mathcal{O}}
\def\cF{\mathcal{F}}
\def\cP{\mathcal{P}}
\def\cV{\mathcal{V}}
\def\cL{\mathcal{L}}
\def\cW{\mathcal{W}}

\def\fm{\mathfrak{m}}

\def\clap#1{\hbox to 0pt{\hss#1\hss}}

\def\mclap{\mathpalette\mathclapinternal}

\def\mathclapinternal#1#2{%
\clap{$\mathsurround=0pt#1{#2}$}}

\def\ce{\mathrel{\mathop:}=}  

\def\fto{\longrightarrow}
\def\injto{\lhook\joinrel\relbar\!\!\:\joinrel\rightarrow}
\def\surjto{\relbar\joinrel\twoheadrightarrow}

\DeclareMathOperator{\Div}{Div}
\DeclareMathOperator{\CDiv}{CDiv}
\DeclareMathOperator{\Pic}{Pic}
\DeclareMathOperator{\Cl}{Cl}
\DeclareMathOperator{\Neg}{neg}
\DeclareMathOperator{\link}{link}
\DeclareMathOperator{\im}{im}
\DeclareMathOperator{\lcm}{lcm}
\DeclareMathOperator{\sign}{sign}

\newcommand{\mailurl}[1]{\href{mailto:#1}{#1}}

\begin{document}

\baselineskip=14pt
\parskip 5pt plus 1pt 

\title{Computing Cohomology on Toric Varieties}
\author{Benjamin Jurke}
\address{Northeastern University, Department of Physics, Dana Research Center, 110 Forsyth Street, Boston, MA 02115, USA}
\address{Max-Planck-Institut f\"ur Physik, F\"ohringer Ring 6, 80805 M\"unchen, Germany}
\email{\mailurl{b.jurke@neu.edu}, \mailurl{mail@benjaminjurke.net}}
\urladdr{\url{http://benjaminjurke.net}}
\date{\today}

\subjclass[2000]{14M25 (Primary); 13D45, 14Q99 (Secondary)}
\keywords{Cohomology of line bundle; Toric variety; Alexander duality}
\thanks{{\it Report number.} MPP-2011-106}


\begin{abstract}
In these notes a recently developed technique for the computation of line bundle-valued sheaf cohomology group dimensions on toric varieties is reviewed. The key result is a vanishing theorem for the contributing components which depends on the structure of the Stanley-Reisner ideal generators. A particular focus is placed on the (simplicial) Alexander duality that provides a central tool for the two known proofs of the algorithm.
\end{abstract}

\maketitle


\section{Introduction \& Motivation}
Cohomology groups play a central role in string model building, where they determine numerous critical properties like the (chiral) zero mode spectrum, Yukawa couplings or counting the number of moduli, for example. The availability of efficient methods to deal with this computational problem is therefore an important requirement.

The majority of geometries considered in theoretical physics and (string) model building is based on toric geometry \cite{Fulton, CoxLittleSchenk, Kreuzer:2006ax}. Due to the requirement of $\cN=1$ supersymmetry in the effective fourdimensional theory, one encounters compact Calabi-Yau threefolds and fourfolds as the compactification spaces. These are typically constructed as hypersurfaces or complete intersections of hypersurfaces in toric varieties. An additional ingredient in string model building are background fluxes, which are described by vector bundles over the compactification space. Such vector bundles can be constructed via different methods, but for computational practicality three types are distinguished:
\begin{itemize}
  \item the \emph{monad bundle} construction, which using a short exact sequence constructs a non-trivial vector bundle $\cV$ on the toric variety $X$ from two other bundles that are typically chosen to be Whitney sums of line bundles:
  \beq
    0 \fto \bigoplus_{i=1}^n \cO_X(a_i) \stackrel{f}{\injto} \bigoplus_{j=1}^m \cO_X(b_j) \stackrel{g}{\surjto} \cV \fto 0
  \eeq
  \item the \emph{extension bundle} construction, which is rather similar but in practice computationally often much more difficult to handle:
  \beq
    0 \fto \bigoplus_{i=1}^n \cO_X(a_i) \stackrel{f}{\injto} \cW \stackrel{g}{\surjto} \bigoplus_{k=1}^l \cO_X(c_k) \fto 0
  \eeq
  \item the \emph{spectral cover} construction, which produces stable holomorphic vector bundles of $SU(n)$ structure group on elliptically-fibered Calabi-Yau threefolds via the Fourier-Mukai transformation from line bundles on top of the so-called spectral cover.
\end{itemize}
In all three of those constructions line bundles provide the basic building block, thus the line bundle-valued cohomology is necessarily involved. In the setting of toric geometry one can relate the line bundles on a hypersurface (or complete intersection of hypersurfaces) to the line bundles of the ambient space via the Koszul sequence
\beq
  0 \fto \cO_X(-S)\injto \cO_X \surjto \cO_S \fto 0.
\eeq
Since short exact sequences of bundles induce long exact sequences of cohomology groups, in the end the critical starting point of every cohomology computation is the knowledge of line bundle-valued cohomology groups on the ambient toric variety.

In \cite{Blumenhagen:2010pv} a novel technique for the computation of the line bundle-valued cohomology group dimensions $h^i(X;\cL_X)$ was introduced and later rigorously proven in \cite{2010arXiv1006.0780J, Rahn:2010fm}. Applications like the constructions above and generalizations surpassing the original scope have been in detail discussed in \cite{Blumenhagen:2010ed} and are summarized in \cite{Blumenhagen:2011xn}. A high-performance implementation called \cohomCalg \cite{cohomCalg:Implementation} was provided along with the original conjecture of the algorithm and has subsequently been improved and optimized. Note that various alternative approaches have been known for some time, like e.g.~\S3.5 of \cite{Fulton}, \cite{EisenbudMustataStillman}, prop.~4.1 in \cite{Borisov} or \S9.1 of \cite{CoxLittleSchenk}.

\subsubsection*{Summary of contents}
It is the goal of these notes to constructively follow the basic structures of the proofs \cite{2010arXiv1006.0780J, Rahn:2010fm} and highlight certain mathematical properties of the algorithm. In section~\ref{sec:ToricGeomAlexDual} the basic notions of toric geometry and the Alexander duality are introduced. Section~\ref{sec:LocalCohomAndGrading} summarizes the correspondence between sheaves and modules and shows how graded components of local cohomology can be identified with line bundle-valued sheaf cohomology groups. Section~\ref{sec:MultiplicityFactors} shows how certain graded components can be grouped together and defines multiplicity factors. In section~\ref{sec:SRidealReductions} the key vanishing result of the algorithm is formulated and section~\ref{sec:MultiplicityComputation} explains an efficient method to compute the multiplicity factors from certain simplicial complexes while highlighting some subleties between the two proofs of the conjecture.

\section{Simplicial Alexander Duality and Toric Geometry}\label{sec:ToricGeomAlexDual}
Consider a finite vertex set $V$ and an \emph{(abstract) simplicial complex} $\Delta$, which is a set of subsets of $V$ such that for each $\sigma\in\Delta$ the subsets $\tau\subset\sigma$ are contained as well, i.e.~$\tau\in\Delta$. In the context of toric geometry $\Delta$ can be identified with the toric fan $\Sigma$, where the $\sigma\in\Delta\cong\Sigma$ are called \emph{cones} and the statement is then simply that all faces of a cone are themself cones of the simplicial complex. Given a set of vertices $\sigma\in V$ let $\hat\sigma\ce V\setminus\sigma$ denote the complement vertices. For each $\sigma\subset V$ the \emph{restriction} of a simplicial complex $\Delta$ on $V$ is defined by $\Delta|_\sigma\ce\{\tau\in\Delta:\tau\subseteq\sigma\}$. 

The \emph{Alexander dual simplicial complex} of $\Delta$ on $V$ is then defined by
\beq\label{eq:SimplicialAlexanderDual}
  \Delta^*\ce\{\sigma\subseteq V : \hat\sigma\not\in\Delta\},
\eeq
i.e.~it consists of all sets of vertices whose complement are not cones of the original simplex. Note that $\Delta^*$ itself defines an (abstract) simplicial complex on $V$. The Alexander dual is a true duality in the sense that $(\Delta^*)^*=\Delta$ and there are examples of self-dual simplices, where $\Delta^*=\Delta$ after a vertex relabeling. The \emph{simplicial Alexander duality} then provides that for each $i$ there exists an isomorphism such that
\beq\label{eq:RestrictedSimplicialAlexanderDual}
  \tilde H_i(\Delta^*) \cong \tilde H^{|V|-3-i}(\Delta),
\eeq
identifying the (reduced) simplicial homology of the complex with the cohomology of its Alexander dual, which can be treated as a standard combination of Poincar\'{e} duality and excision. Several detailed examples can be found in e.g.~chap.~5 of \cite{MillerSturmfels} and further information in \cite{Moller:AlexDuality}.

Given a cone $\sigma\in\Delta$, we can also define the \emph{link of $\sigma$ in $\Delta$} by
\beq
  \link_\Delta(\sigma) \ce \{ \tau \in\Delta : \tau\cup\sigma\in\Delta, \; \tau\cap\sigma=\emptyset\}.
\eeq
Note that this defines an actual simplicial complex on $\Delta|_{\hat\sigma}$. The simplicial Alexander duality \eqref{eq:SimplicialAlexanderDual} can then be restated as
\beq\label{eq:RestrictedAlexanderDuality}
  \tilde H_i\big(\link_{\Delta^*}(\sigma)\big) \cong \tilde H^{|V|-|\sigma|-3-i}(\Delta|_{\hat\sigma}).
\eeq

The Alexander dual can also be formulated in terms of ideals and squarefree monomials, which in fact offers a more convenient perspective in the context of toric geometry. We consider a $d$-dimensional simplicial projective toric variety $X$ and a fan $\Sigma$ in the lattice $N\cong\IZ^d$. As mentioned before, one can treat a fan $\Sigma$ as a simplicial complex $\Delta$ and we will do so from this point on. Let $x_1,\dots,x_n$ with $n=|V|$ be homogeneous coordinates associated to the vertices in $V$ that generate the rays (1d cones) of the fan and let $\boldsymbol{x}^\sigma\ce\prod_{i\in\sigma} x_i$ denote the associated squarefree monomial for some $\sigma\subseteq{}[n]\ce\{1,\dots,n\}$ after labeling the vertices via $[n]{}\cong V$. Then let
\beq
  \fm^\sigma\ce\langle x_i : i\in\sigma \rangle
\eeq
be the monomial prime ideal corresponding to the vertices $\sigma$ in the so-called \emph{Cox ring} $S\ce\IC[\boldsymbol{x}]=\IC[x_1,\dots,x_n]$ of homogeneous coordinates. Given a monomial ideal $J=\langle \boldsymbol{x}^{\sigma_1},\dots,\boldsymbol{x}^{\sigma_r}\rangle$ in $S$, the \emph{Alexander dual monomial ideal} is
\beq
  J^* = \fm^{\sigma_1}\cap\dots\cap\fm^{\sigma_r}.
\eeq

As before let $\hat\sigma\ce{}[n]{}\setminus\sigma$ denote the complement of some $\sigma\subseteq{}[n]{}\cong V$. Then
\beq
  B_\Sigma \ce \langle \boldsymbol{x}^\sigma : \hat\sigma\in\Sigma \rangle
\eeq
defines the \emph{irrelevant ideal} of the fan $\Sigma$, and a minimal generating set for $B_\Sigma$ is given by the monomials corresponding to the complements of the maximal cones of $\Sigma$. The \emph{Stanley-Reisner ideal} on the other hand is defined by
\beq\label{eq:StanleyReisnerDef}
  I_\Sigma\ce\langle \boldsymbol{x}^\sigma : \sigma\not\in\Sigma \rangle,
\eeq
which---using appropiate identifications---is obviously Alexander dual to the irrelevant ideal by \eqref{eq:SimplicialAlexanderDual}:
\beq\label{eq:IrrVsSRidealAlexDual}
  (B_\Sigma)^* \cong I_\Sigma,
	\qquad
	(I_\Sigma)^* \cong B_\Sigma.
\eeq
The Stanley-Reisner ideal and indirectly the associated \emph{Stanley-Reisner ring} $S/I_\Sigma$ take a central role in our computational technique. Likewise, one can also consider the ring $S/B_\Sigma$ which in an abuse of the terminology can be treated as the Alexander dual to $S/I_\Sigma$ based on \eqref{eq:IrrVsSRidealAlexDual}. See \cite{Miller:CohenMacaulayCriteria} for further information on the various relationships between (squarefree) monomial ideals and simplicial complexes.

\section{Local Cohomology and Grading Reorganization}\label{sec:LocalCohomAndGrading}
Consider the Weil divisor class group $\Cl(X)=\Div(X)/\Div_0(X)$ of $X$. Since we are considering a smooth variety $X$ all Weil divisors are also Cartier,\footnote{\emph{Weil divisors} $\Div(X)$ are formal sums of irreducible codimension-1 subvarieties of $X$. A Weil divisor is called a \emph{Cartier divisor} $\CDiv(X)$ if it is locally principal, i.e.~when it can be locally described by the vanishing order and locus of a rational function.} such that the class group can be identified with the Picard group $\Pic(X)=\CDiv(X)/\Div_0(X)$ and therefore $\Cl(X)\cong\Pic(X)\cong\IZ^{n-d}$.

The Cox ring $S=\IC[\boldsymbol{x}]=\IC[x_1,\dots,x_n]$ of homogeneous coordinates is a $\Cl(X)$-graded ring, i.e.~it can be decomposed like
\beq
  S=\bigoplus_{\mclap{\alpha\in\Cl(X)}}S_\alpha, \qquad \text{such that $S_\alpha \cdot S_\beta \subset S_{\alpha+\beta}$},
\eeq
and each of those graded decomposition spaces is naturally isomorphic to the space of sections of the line bundle $\cO_X(\alpha)$, i.e.
\beq
  S_\alpha \cong \Gamma\big(X; \cO_X(\alpha)\big).
\eeq
This identification forms the basis of a deeper connection between line bundle-valued cohomology and algebraic notions. Given an $S$-module $M$ and an ideal $J\subset S$ the \emph{$J$-torsion submodule} is defined by
\beq
  \Gamma_J(M) \ce \{ \boldsymbol{x}\in M : J^k \cdot \boldsymbol{x} = 0 \text{ for some $k\in\IN$}\}.
\eeq
The \emph{$i$-th local cohomology} $H^i_J(M)$ of $M$ with support on $J$ is then defined by the $i$-th cohomology of the complex
\beq
  0 \fto \Gamma_J(I^0) \fto \Gamma_J(I^1) \fto \Gamma_J(I^2) \fto \ldots
\eeq
that is obtained from an injective resolution $0\fto I^0\ce M \fto I^1 \fto \ldots$ of the module $M$, see \S9.5 of \cite{CoxLittleSchenk} for a detailed introduction of the $\Gamma_J(\,\cdot\,)$ functor.

The local cohomology also inherits any grading of $M$ and the precise connection between line bundle cohomology and local cohomology can then be formulated as
\beq\label{eq:CoarseGrading}
  H^i\big(X;\cO_X(\alpha)\big)\cong H^{i+1}_{B_\Sigma}(S)_\alpha
\eeq
for any divisor class $\alpha\in\Cl(X)\cong\Pic(X)$ and $i\ge1$,\footnote{There seems to be a small error in \cite{Rahn:2010fm} in the argumentation surrounding eqn.~10, but since the computation of $h^0\big(X;\cO_X(\alpha)\big)=\dim_\IC\Gamma\big(X;\cO_X(\alpha)\big)$ is equivalent to counting only $|(\alpha,\emptyset)|$ with multiplicity factor 1 (see \eqref{eq:CountingNums} and \eqref{eq:AlgorithmFormula} below) the end result remains unchanged.} see prop.~2.3 of \cite{EisenbudMustataStillman} or thm.~9.5.7 in \cite{CoxLittleSchenk} for a proof. Computing the $i$-th line bundle-valued cohomology group of a toric variety therefore is equivalent to computing the $(i+1)$-th cohomology group of the homogeneous coordinate ring $S$ localized on the irrelevant ideal $B_\Sigma$, and the line bundle $\cO_X(\alpha)$ determines which graded piece has to be considered.

The $\Cl(X)$-grading of $S$ and $H^i_{B_\Sigma}(S)$ can be refined by introducing a $\IZ^n$-grading induced from projective weights of the homogeneous coordinates $x_i$ themselves, which in the physics literature are often referred to as the GLSM charges $Q_i^{(j)}$. More precisely, we consider the map 
\beq
  \bal
	  f:{}&{}\IZ^n \fto \Cl(X)\cong\IZ^{n-d} \\
		&{} \vec e_i \mapsto [D_i] = (Q_i^{(1)},\dots,Q_i^{(n-d)}),
	\eal
\eeq
where $\vec e_i\in\IZ^n$ is a basis vector associated to the coordinate $x_i$, such that a monomial $x_1^{k_1}\cdots x_n^{k_n}$ can be simply represented as $k_1\vec e_1+\ldots+ k_n\vec e_n$. The coordinate divisor $D_i$ refers to the hypersurface $\{x_i=0\}\subset X$ and $[D_i]\in\Cl(X)\cong\Pic(X)$ to its divisor class. In terms of this finer grading we can then use
\beq\label{eq:FineSplitting}
  H^i\big(X;\cO_X(\alpha)\big) = \bigoplus_{\mclap{\substack{\vec u\in\IZ^n: \\ f(\vec u)=\alpha}}} H^{i+1}_{B_\Sigma}(S)_{\vec u},
\eeq
which means that we split up $H^{i+1}_{B_\Sigma}(S)_\alpha$ from \eqref{eq:CoarseGrading} into $f^{-1}(\alpha)\subset\IZ^n$ pieces. In other words: For the grading we are considering all monomials $x_1^{k_1}\cdots x_n^{k_n}$ whose total projective weight/GLSM charge/degree is equal to $\alpha\in\Cl(X)$ that specifies the line bundle $\cO_X(\alpha)$. An explicit way to compute this local cohomology $H^i_{B_\Sigma}(S)$ directly via generalized \v{C}ech cohomology of a free resolution of the \emph{irrelevant ring} $S/B_\Sigma$ is summarized in \S3 of \cite{2010arXiv1006.0780J}, but we will only use it as an intermediate step---however, see the comments at the end of sec.~\ref{sec:MultiplicityComputation} below.

It remains to understand the structure of the splitting \eqref{eq:FineSplitting} better and ideally to group the non-vanishing summands together. The graded pieces $H^{i+1}_{B_\Sigma}(S)_{\vec u}$ can also be computed in terms of an abstract simplicial complex on the vertex set/coordinates of the fan: For $\vec u\in\IZ^n$ let
\beq
  \Neg(\vec u) \ce \{ k \in V \cong {}[n]{} : u_k < 0 \}
\eeq
be the set of those indices where $\vec u$ has a negative entry and as usual define the complement by $\widehat{\Neg}(\vec u) \ce {} [n]\setminus \Neg(\vec u) = \{ k \in {}[n] : u_k \ge 0\} \subseteq {}[n]$. Following prop.~3.1 in \cite{2010arXiv1006.0780J} it can then be shown that
\beq\label{eq:NegIdentificationIso}
  H^{i+1}_{B_\Sigma}(S)_{\vec u} \cong \tilde H_{d-i-1}(\Sigma|_{\widehat{\Neg}(\vec u)})
\eeq
which in particular implies that $H^i_{B_\Sigma}(S)_{\vec u}\cong H^i_{B_\Sigma}(S)_{\vec v}$ if $\Neg(\vec u)=\Neg(\vec v)$. The refined $\IZ^n$-grading in \eqref{eq:FineSplitting} can therefore be simplified by collecting the pieces $H^i_{B_\Sigma}(S)_{\vec u}$ with the same $\Neg(\vec u)$. 

Following the notation introduced in \cite{Rahn:2010fm}, given some $\sigma\subseteq{}[n]\cong V$ define
\beq\label{eq:AssociatedZnGrading}
  \tilde\sigma\in\IZ^n \quad \text{such that} \quad \tilde\sigma_i \ce \begin{cases} 1 & \text{if $i\in\sigma$,} \\ 0 & \text{if $i\not\in\sigma$,} \end{cases}
\eeq
i.e.~we have a tupel of 0s or 1s depending on whether or not the respective coordinate index appears in $\sigma$. For the computation of the dimension of $H^i\big(X;\cO_X(\alpha)\big)$ based on \eqref{eq:FineSplitting} and using $\smash{H^i_{B_\Sigma}(S)_{\vec u}\cong H^i_{B_\Sigma}(S)_{-\widetilde{\Neg}(\vec u)}}$ we therefore arrive at
\beq\label{eq:CollectedFormula}
  h^i\big(X;\cO_X(\alpha)\big) = \sum_{\mclap{\sigma\subseteq{}[n]}} |(\alpha,\sigma)|\cdot h^{i+1}_{B_\Sigma}(S)_{-\tilde\sigma}
\eeq
due to \eqref{eq:NegIdentificationIso}, where $|(\alpha,\sigma)|$ counts the number of identical $\IZ^n$-pieces $H^i_{B_\Sigma}(S)_{\vec u}$ via
\beq\label{eq:CountingNums}
  (\alpha,\sigma) \ce \{ \vec u \in \IZ^n : f(\vec u)=\alpha,\;\Neg(\vec u)=\sigma \}.
\eeq
The potentially infinite sum in \eqref{eq:FineSplitting} has therefore been reduced to a finite sum of $2^n$ terms, provided that the counting of \eqref{eq:CountingNums} can be handled efficiently.\footnote{In the original conjecture \cite{Blumenhagen:2010pv} the factor $|(\alpha,\sigma)|$ appears as the counting of ``rationoms''/rational functions which refers to the Laurent monomials. At this point it remains to understand the $\smash{h^{i+1}_{B_\Sigma}(S)_{-\tilde\sigma}}$, which will turn out to be the secondary/remnant cohomology factors along with a further reduction in computational complexity.}

\section{Graded Betti Numbers, Resolutions and Multiplicity Factors}\label{sec:MultiplicityFactors}
In order to restrict the computational effort in \eqref{eq:CollectedFormula} further another ingredient is required. We consider the minimal free resolution of a finitely generated graded $S$-module $M$, which can be written in the form
\beq
  \cF_\bullet:\quad 0 \longleftarrow V_0 \stackrel{\phi_1}{\longleftarrow} V_1 \longleftarrow \dots \longleftarrow V_{\ell-1} \stackrel{\phi_\ell}{\longleftarrow} V_\ell \longleftarrow 0
\eeq
and which is \emph{acyclic}, meaning that it is exact everywhere except in the first position where $M \cong V_0 / \im(\phi_1)$ holds. With respect to a $\IZ^n$-grading the individual spaces $V_i$ of the resolution can be written as
\beq
  V_i = \bigoplus_{\mclap{{\vec u}\in\IZ^n}} \big(S_{-\vec u}\big)^{\beta_{i,\vec u}},
\eeq
which also defines the \emph{graded Betti numbers} $\beta_{i,\vec u}$, see e.g.~\cite{bruns1993cohen-macaulay} for details. A \emph{minimal} resolution minimizes the ranks of the graded $S$-modules $V_i$ and is unique up to isomorphisms. Furthermore, any $S$-module always has a free resolution with length $\ell$---assuming that $V_\ell\not=0$---smaller or equal $n$.

By treating both the Stanley-Reisner ideal $I_\Sigma$ and the Stanley-Reisner ring $S/I_\Sigma$ as $\IZ^n$-graded $S$-modules, one can consider their minimal free resolutions. The associated graded Betti numbers can be calculated by the \emph{Hochster formula}
\beq
  \bal
    \beta_{i-1,\tilde\sigma}(I_\Sigma) {}&{} = \beta_{i,\tilde\sigma}(S/I_\Sigma) \\
		{}&{}= \dim \tilde H^{|\sigma|-i-1}(\Sigma|_\sigma) \;\smash{\stackrel{\eqref{eq:RestrictedAlexanderDuality}}{=}}\; \dim \tilde H_{i-2}\big(\link_{\Delta^*}(\hat\sigma)\big),
	\eal
\eeq
see e.g.~cor.~5.12 in \cite{MillerSturmfels}, which allows to rewrite \eqref{eq:CollectedFormula} in the form
\beq\label{eq:UnSRreducedFormula}
  h^i\big(X;\cO_X(\alpha)\big) = \sum_{\mclap{\sigma\subseteq{}[n]}} |(\alpha,\sigma)|\cdot\beta_{|\sigma|-i,\tilde\sigma}(S/I_\Sigma).
\eeq
However, we still have to evaluate all possible subsets $\sigma\subseteq{}[n]\cong V$, i.e.~$2^n$ terms, to compute $h^i\big(X;\cO_X(\alpha)\big)$.

\section{Reductions based on Stanley-Reisner Ideal Generators}\label{sec:SRidealReductions}
The central point of our algorithm \cite{Blumenhagen:2010pv} concerns a further reduction of the number of terms in \eqref{eq:UnSRreducedFormula} by taking the Stanley-Reisner ideal $I_\Sigma$ into account. More precisely, following definition \eqref{eq:StanleyReisnerDef} let $I_\Sigma=\langle S_1,\dots,S_t\rangle$ be generated by $t$ squarefree monomials $S_r = \boldsymbol{x}^{\sigma_r}$ in the Cox ring $S$. We basically need to consider all possible subsets of $I_\Sigma$-generators and consider the union of coordinates appearing in the generators. In more formal terms: for all $\tau\subseteq{}[t]$ denote the associated $\IZ^n$-degree by
\beq\label{eq:SRunionZnDegreeDef}
  \vec a_\tau \ce \deg_{\boldsymbol x}(S_\tau)\in\IZ^n
	\quad\text{where}\quad
	S_\tau \ce \lcm_{\boldsymbol{x}}\{S_i : i\in\tau\}
\eeq
represents this ``union of coordinates''---the least common multiple of several generators with respect to the homogeneous coordinates $x_1,\dots,x_n$ of the Cox ring. Note that $\vec a_\tau$ is also a series of 0s and 1s like \eqref{eq:AssociatedZnGrading}.

Based on this we consider the set of all such $\vec a_\tau$, basically the set of all $\IZ_n$-degrees that can arise from unions of Stanley-Reisner ideal generators $S_i$, which will be denoted by
\beq
  \cP(I_\Sigma) \ce \{ \vec a_\tau : \tau\subseteq{}[t] \}.
\eeq
Note that different $\tau\subseteq{}[t]$ may lead to the same $\IZ^n$-degree $\vec a_\tau$. Coming back to the computation \eqref{eq:UnSRreducedFormula}, the most important result is that we only have to consider combinations of the coordinates $\sigma\subseteq{}[n]\cong V$ that are found in $\cP(I_\Sigma)$:
\begin{itemize}
  \item All collections of coordinates that are not unions of coordinates of Stanley-Reisner ideal generators do not contribute to $h^i\big(X;\cO_X(\alpha)\big)$, i.e.
	      \beq
				  \boxed{\beta_{r,\tilde\sigma}(S/I_\Sigma) = 0 \quad\text{for all}\quad \sigma\subseteq{}[n] \text{ where } \tilde\sigma\not\in\cP(I_\Sigma).}
				\eeq
\end{itemize}
This leads to the original algorithm formula of \cite{Blumenhagen:2010pv}. Due to a further ''Serre duality for graded Betti numbers``, that was first observed in \cite{Rahn:2010fm} and proven in \cite{2010arXiv1006.0780J}, the same is true for the complement $\hat\sigma$ such that the sum in \eqref{eq:UnSRreducedFormula} is in the end reduced to
\beq\label{eq:AlgorithmFormula}
  \boxed{
  h^i\big(X;\cO_X(\alpha)\big) = \sum_{\mclap{\substack{\sigma\subseteq{}[n] \\ \tilde\sigma,\tilde{\hat\sigma}\in\cP(I_\Sigma)}}} \overbrace{|(\alpha,\sigma)|}^{\mclap{\text{elements in ``neg-group''}}}\cdot\underbrace{\beta_{|\sigma|-i,\tilde\sigma}(S/I_\Sigma)}_{\mclap{\text{multiplicity factor}}}}
\eeq
with $(\alpha,\sigma)$ defined as in \eqref{eq:CountingNums}. Formally \eqref{eq:AlgorithmFormula} always represents a reduction compared to \eqref{eq:UnSRreducedFormula} that strongly depends on the form of the Stanley-Reisner ideal $I_\Sigma$.

\section{Computing Multiplicity Factors}\label{sec:MultiplicityComputation}
While the counting of the number of neg-group elements $|(\alpha,\sigma)|$ is the well-known task of solving a linear system over the integers, it remains to provide an efficient method for the computation of the graded Betti numbers $\beta_{|\sigma|-i,\tilde\sigma}(S/I_\Sigma)$ which have been dubbed multiplicity factors in \cite{Blumenhagen:2010pv, 2010arXiv1006.0780J, Rahn:2010fm, Blumenhagen:2010ed, Blumenhagen:2011xn}. This can be done by defining a subcomplex of the full abstract simplicial complex $\Delta_{[t]}$ by taking only those $\tau\subseteq{}[t]$ whose $\IZ^n$-degree $\vec a_\tau$ has a certain value. More precisely, for a given $\sigma\subseteq{}[n]$ define a simplicial (sub)complex by
\beq\label{eq:TaylorSubcomplex}
  \Gamma^\sigma \ce \{ \tau\subseteq{}[t] : \vec a_\tau = \tilde\sigma \}.
\eeq
As mentioned before, different $\tau\subseteq{}[t]$ share the same $\IZ^n$-degree $\vec a_\tau$ according to \eqref{eq:SRunionZnDegreeDef}, such that the number of elements in $\tau$---the cone dimension---becomes a major distinguishing aspect. Let $F_j(\Gamma^\sigma)$ denote the \emph{$j$-faces of $\Gamma^\sigma$}, i.e.~the $\IC$-vector space whose basis vectors $e_\tau$ are labeled by the $\tau\in\Gamma^\sigma$ having $|\tau|=j+1$ elements. The complex mappings correspond to the standard boundary mappings in a simplicial complex, specifically by linearly extending
\beq
  \bal
	  \phi_j : {} & F_j(\Gamma^\sigma) \fto F_{j-1}(\Gamma^\sigma) \\
		&{} e_\tau \mapsto \sum_{i\in\tau} \sign(i,\tau)\; e_{\tau\setminus\{i\}},
	\eal
\eeq
where $\sign(i,\tau) \ce (-1)^{s-1}$ when $i$ is the $s$-th element of $\tau\subset{}[t]=\{1,\dots,t\}$ written in increasing order. The second central proposal of \cite{Blumenhagen:2010pv} can then be stated as follows:
\begin{itemize}
  \item The graded Betti numbers $\beta_{r,\tilde\sigma}(S/I_\Sigma)$ appearing in \eqref{eq:AlgorithmFormula} can be computed from the reduced homology dimension of the simplicial complex $\Gamma^\sigma$ defined in \eqref{eq:TaylorSubcomplex}, i.e.
	      \beq
				  \boxed{\beta_{r,\tilde\sigma}(S/I_\Sigma) = \beta_{r-1,\tilde\sigma}(I_\Sigma) = \dim_\IC \tilde H_{r-1}(\Gamma^\sigma).}
				\eeq
\end{itemize}
Due to the basic structure of the simplicial complex $\Gamma^\sigma$ the evaluation of this reduced homology is straightforward, thus providing a convenient method to determine the multiplicity factors required in \eqref{eq:AlgorithmFormula}.

In \cite{2010arXiv1006.0780J} an analogous computation of the graded Betti numbers/multiplicity factors is carried out using the abstract simplicial complex $\Lambda_I$, which is complementary to $\Gamma^\sigma$ within the full simplex $\Delta_{[m]}$ of all $\boldsymbol{x}^\sigma$-dividing $I_\Sigma$-generators $S_1,\dots,S_m$, see eqn.~46 in \cite{Rahn:2010fm} for the precise relationship. Moreover, the approach in \cite{2010arXiv1006.0780J} is based on computing the local cohomology $\smash{H^i_{B_\Sigma}(S)}$ via a free resolution of $S/B_\Sigma$, whereas the proof in \cite{Rahn:2010fm} is based on a free resolution of $S/I_\Sigma$. From this perspective both proofs \cite{Rahn:2010fm} and \cite{2010arXiv1006.0780J} can be considered to be complementary: while \cite{2010arXiv1006.0780J} is somewhat more concise, the proof in \cite{Rahn:2010fm} by considering the complex $\Gamma^\sigma$ more closely follows the original conjecture in \cite{Blumenhagen:2010pv}. The key relation between the two approaches is obviously the Alexander duality according to \eqref{eq:IrrVsSRidealAlexDual}.\footnote{It should be mentioned that the generalized \v{C}ech complex associated to the resolution of the (Alexander dual) ring $S/B_\Sigma$ in \cite{2010arXiv1006.0780J} is closely related to the toric \v{C}ech complex in eqn.~4.1 of \cite{Borisov}: A toric variety can be patched together from charts on
\beq
  U_\sigma \ce \big\{ (x_1,\dots,x_n) : \text{$x_i = 0$ only for $i\in\sigma$} \big\},
\eeq
i.e.~on $U_\sigma$ one can consider the ring of Laurent monomials $S\big[ \frac{1}{\boldsymbol{x}^{\hat\sigma}} \big]$ which allows for negative exponents in those coordinates $x_i$ where $i\not\in\sigma$. This coincides with the space of sections of the holomorphic line bundle $\cO$ on $U_\sigma$.}

Note that in actual computations one can determine both $\cP(I_\Sigma)$ and the relevant simplicial complexes $\Gamma^\sigma$ in a single step. However, since it is necessary to evaluate the entire powerset of Stanley-Reisner ideal generators with its $2^t$ elements, the algorithm complexity---aside from the counting of $|(\alpha,\sigma)|$ and actually computing the reduced cohomology of $\Gamma^\sigma$---grows exponentially with the number of Stanley-Reisner ideal generators. This behaviour is also seen in the reference implementation {\cohomCalg} \cite{cohomCalg:Implementation}. Compared to the starting point \eqref{eq:CollectedFormula} one effectively exchanges an exponential growth in the number of vertices $n$ by an exponential growth in the number of Stanley-Reisner ideal generators $t$.

When the number of Stanley-Reisner ideal generators becomes large, the algorithm and its implementation reach their limits. Other tools capable of computing line bundle-valued sheaf cohomology group dimensions include the {\sl Sage} system \cite{SAGE} or the {\sl Macaulay2} \cite{M2} packages ``{\tt NormalToricVarities}'' by G.~Smith or ``{\tt ToricVectorBundles}'' by R.~Birkner, N.~O.~Ilten and L.~Petersen.

\subsection*{Acknowledgments}
The author would like to thank L.~Borisov for some helpful discussions regarding the importance of the Alexander duality in the algorithm and would like to acknowledge the University of Pennsylvania for hospitality during the String Math 2011 conference. The algorithm discussed herein is based on the collaborative efforts together with R.~Blumenhagen, T.~Rahn and H.~Roschy and the author would like to thank T.~Rahn for discussions on the manuscript. This work was partially supported by the NFS-Microsoft grant NSF/CCF-1048082.


\bibliography{rev1}
\bibliographystyle{utphys}

\end{document}